\documentclass{elsarticle}   	
\usepackage{geometry}                		
\usepackage[parfill]{parskip}    		
\usepackage{graphicx}				
\usepackage{amsmath,amsfonts,amssymb}
\usepackage{color}

\DeclareSymbolFont{bbold}{U}{bbold}{m}{n}
\DeclareSymbolFontAlphabet{\mathbbold}{bbold}

\newcommand{\B}[1]{\boldsymbol{#1}}
\newcommand{\prt}[2]{{#1_{#2}}}

\newcommand{\R}{\mathbf{R}}

\newcommand{\U}{\mathcal{U}}
\newcommand{\T}{\mathcal{T}}
\newcommand{\Ell}{\mathcal{L}}

\newcommand{\Du}{\delta_u}
\newcommand{\Dbar}{\bar{\delta}_u}
\newcommand{\FD}{\delta_t}

\newcommand{\dt}{\tau}

\begin{document}
\begin{frontmatter}

\title{Symplectic Runge-Kutta discretization of a regularized forward-backward sweep iteration for optimal control problems}
\author{Xin Liu\corref{cor1}\fnref{csc}}
\cortext[cor1]{Corresponding author.}
\ead{x.liu2@uu.nl}
\author{Jason Frank\corref{}}
\ead{j.e.frank@uu.nl}
\address{Mathematical Institute, Utrecht University, P.O. Box 80010, 3508 TA Utrecht, the Netherlands}
\fntext[csc]{The first author gratefully acknowledges support from the Chinese Scholarship Council under grant number 201607040074.}

\begin{abstract}
Li, Chen, Tai \& E. (\emph{J.~Machine Learning Research}, 2018)  have proposed a regularization of the forward-backward sweep iteration for solving the Pontryagin maximum principle in optimal control problems. The authors prove the global convergence of the iteration in the continuous time case.  In this article we show that their proof can be extended to the case of numerical discretization by symplectic Runge-Kutta pairs. We demonstrate the convergence with a simple numerical experiment.
\end{abstract}

\end{frontmatter}

Recently, Li et al.~\cite{LiChTaE18} proposed a new indirect iteration for optimal control problems in the context of deep neural networks, that utilizes the `method of successive approximations', i.e.~forward and backward integrations, combined with an `augmented Lagrangian' regularization that ensures global convergence.  The authors argue that this approach is particularly suitable for high-dimensional optimal control problems as encountered in deep learning. Large scale optimal control problems figure centrally in a number of modern applications such as deep neural networks \cite{LiChTaE18}, reinforcement learning \cite{SuBa18,Bersekas}, filtering and data assimilation methods \cite{BaCr09,ZhFrHeSh15} and mean field and stochastic differential games \cite{CaDe18}.
In this paper we describe how the iteration of Li et al.~combines naturally with symplectic/variational integrators to yield a convergent numerical scheme.

Optimal control problems possess a natural variational structure that gives rise to Hamiltonian dynamics which may be exploited in a numerical treatment \cite{JuMaOb05}. Symplectic methods for Hamiltonian \emph{initial} value problems have been much studied since the mid-1990s due to their demonstrated superiority for conserving energy and other first integrals \cite{SaSeCa94,HaLuWa06,LeRe05}. In contrast, optimal control problems lead to \emph{boundary} value problems, and it is unclear that the advantages of symplectic integrators for IVPs should translate to the BVP setting. Recent papers that address the use of symplectic Runge-Kutta methods for optimal control stress the conservation of quadratic invariants \cite{SanzSerna16,FrZh14} and the persistence of critical orbits in modified equation expansions \cite{ChHaVi09}.  See also recent work on the preservation of bifurcations under symplectic discretization of boundary value problems \cite{McLOf18}.

In the first three sections of the paper we review the Hamiltonian structure of optimal control problems (\S \ref{sec:background}),
the regularized forward-backward sweep iteration proposed by Li et al.~\cite{LiChTaE18} (\S \ref{sec:msa}) and the discrete variational approach to constructing symplectic Runge-Kutta methods (\S \ref{sec:VarInt}).  In Section \ref{sec:convergence} we prove the convergence of the discrete regularized forward-backward sweep iteration, which follows closely the proof of \cite{LiChTaE18} for the continuous case. It is the symplectic structure of the discretization that facilitates this proof. Finally, in Section \ref{sec:numer} we demonstrate the convergence of the method for a simple example using two symplectic discretizations.

\section{Background \label{sec:background}}
In this section we define continuous optimal control of differential equations and discuss their Hamiltonian structure, and we review the regularized forward-backward sweep iteration of Li et al.~\cite{LiChTaE18}.

\subsection{Hamiltonian structure of optimal control problems}
The state of the system to be controlled is described by a vector $x(t):\T\to\R^d$, where $\T=[0,T]$ represents a time interval.  The control function $u(t)$ is for each $t$ an element of the set of admissable controls $\U\subset\R^m$.
The motion of the system is described by a differential equation
\begin{equation}\label{dynamics}
	\dot{x} (t)= f(x(t),u(t)), \qquad x(0) = \xi,
\end{equation}
where $f:\R^d\times \U \to\R^d$ and $\xi\in\R^d$ is the initial state. The control $u(t)$ is chosen to minimize the objective functional
\begin{equation}\label{cost}
	J[u] = \Phi(x(T)) + \int_0^T h(x(t),u(t)) \, dt,
\end{equation}
where $\Phi:\R^d \to \R$ is the end cost and $h:\R^d\times \U \to\R$ is the running cost. The cost functional \eqref{cost} and the motion \eqref{dynamics} are assumed to be given input to the problem.

In \cite{LiChTaE18} no running cost $h$ is considered. We include it here because it is present in many applications and its treatment is straightforward. As in \cite{LiChTaE18} (cf.~equations (A1) and (A2) of that article) we assume that $\Phi$ and $f$ are twice continuously differentiable with respect to $x$ and satisfy Lipschitz conditions for all $x$,  $x' \in \R^d$, $u\in\U$ and $t\in \T$. We require similar assumptions on $h$:
\begin{equation} \label{contLipschitz} \begin{split}
	&| \Phi(x) - \Phi(x') | + \| \prt{\Phi}{x}(x) - \prt{\Phi}{x}(x') \| \le K \| x - x'\|, \\
	&\| f(x,u) - f(x',u) \| + \| \prt{f}{x}(x,u) -  \prt{f}{x}(x',u) \| \le K \| x - x'\|, \\
	&| h(x,u) - h(x',u) | + \| \prt{h}{x}(x,u) -  \prt{h}{x}(x',u) \| \le K \| x - x'\|,
	\end{split}
\end{equation}
where $\prt{h}{x}$ denotes the vector of partial derivatives of $h$ with respect to $x$ and $\prt{f}{x}$ denotes the Jacobian matrix of partial derivatives of $f$ with respect to $x$. Here and throughout the article, we denote by $\| \cdot \|$ the Euclidean norm on vector spaces.
Note that the solution $x(t)$ of \eqref{dynamics} is well-defined for appropriate $u(t)$ so that we may think of $J$ as a functional essentially depending only on $u(t)$.

The problem can be reformulated as a constrained optimization problem by introducing the Lagrange multiplier function $\lambda(t):\T\to\R^d$ and the Lagrangian functional
\begin{equation}\label{Lagrangian}
	\Ell[x,\lambda,u] = \Phi(x(T)) + \lambda_0^T(x(0)-\xi) + \int_0^T h(x,u) + \lambda^T\left(\dot{x} - f(x,u) \right) \, dt.
\end{equation}
(Throughout the paper we use the transpose and dot product notation interchangeably, whichever is more convenient.)
The variational derivatives of the functional  $\Ell$ with respect to the functions $x(t)$, $\lambda(t)$ and $u(t)$,
denoted $\Ell_x$, $\Ell_\lambda$ and $\Ell_u$, are defined with respect to the $L^2$ inner product.  The first order necessary conditions for an optimum of \eqref{Lagrangian} are given by the Euler-Lagrange equations ($\Ell_x\equiv\Ell_\lambda\equiv\Ell_u\equiv 0$):
\begin{align}
	\dot{x} &= f(x,u), \quad x(0) = \xi, \label{EL_x}\\
	\dot{\lambda} &= -\prt{f}{x} (x,u)^T \lambda + \prt{h}{x}(x,u), \quad \lambda(T) = -\prt{\Phi}{x}(x(T)),\label{EL_lam}\\
	0 &=  \prt{f}{u}(x,u)^T \lambda - \prt{h}{u}(x,u). \label{EL_u}
\end{align}
In particular, if $f$ and $h$ are smooth and $u$ is an optimal control in the interior of $\U$, then it satisfies \eqref{EL_x}--\eqref{EL_u}.
It is convenient to define a function $g(x,\lambda,u)$ for the right side of \eqref{EL_lam}:
\begin{equation}\label{gfun}
	g(x,\lambda,u) = -\prt{f}{x} (x,u)^T \lambda + \prt{h}{x}(x,u).
\end{equation}

A Legendre transform yields the Hamiltonian function
\begin{equation}\label{Hamiltonian}
	H(x,\lambda,u) = \lambda^Tf(x,u) - h(x,u),
\end{equation}
and Hamilton's equations are
\begin{align}
	\dot{x} &= \prt{H}{\lambda}(x,\lambda,u), \label{Ham_x}\\
	\dot{\lambda} &= - \prt{H}{x}(x,\lambda,u), \label{Ham_lam}\\
	0 &= \prt{H}{u}(x,\lambda,u). \label{Ham_u}
\end{align}
Note that minimizing the objective functional $J$ corresponds to maximizing the Hamiltonian with respect to $u$. The condition \eqref{Ham_u} above can be generalized to apply to controls $u(t)$ constrained to lie in $\U$ by replacing  \eqref{Ham_u} with Pontryagin's maximum principle
\begin{align}
	\dot{x} &= f(x,u^*), \quad x(0) = \xi, \label{PMP_x}\\
	\dot{\lambda} &= g(x,\lambda,u^*), \quad \lambda(T) = -\prt{\Phi}{x}(x(T)) \label{PMP_lam}\\	
	u^*(t) &= \arg\max_{u(t)\in\U} H(x,\lambda,u), \quad \forall t\in\T \label{PMP_u}
\end{align}

\subsection{Regularized forward-backward sweep iteration\label{sec:msa}}
Solution of  \eqref{PMP_x}--\eqref{PMP_u} is challenging due to the boundary conditions. One approach is to solve in succession \eqref{PMP_x} for $x(t)$, \eqref{PMP_lam} for $\lambda(t)$ and \eqref{PMP_u} for $u^*(t)$ and iterate. Such a forward-backward sweep iteration typically diverges unless the Lipschitz constant $K$ and the time interval $T$ are small \cite{McMoHa12}.
In a recent article, Li et al.~\cite{LiChTaE18} proposed a modified iteration based on a regularized Lagrangian approach.  They introduce the augmented Hamiltonian function
\begin{equation}\label{augLag}
	\tilde{H}(x,\lambda,u,p,q) = H(x,\lambda,u) - \frac{\rho}{2} \left( \| p - \prt{H}{\lambda} (x,\lambda,u) \|^2 + \| q + \prt{H}{x} (x,\lambda,u) \|^2 \right),
\end{equation}
where $\rho>0$ is a regularization parameter.
Subsequently, the forward-backward sweep iteration is modified to solve consecutively:
\begin{align}
	\dot{x}^{(k+1)} &= \prt{\tilde{H}}{\lambda} (x^{(k+1)},\lambda^{(k)},u^{(k)},\dot{x}^{(k+1)},\dot{\lambda}^{(k)}), \label{MSA_x}\\
	\dot{\lambda}^{(k+1)} &= -\prt{\tilde{H}}{x} (x^{(k+1)},\lambda^{(k+1)},u^{(k)},\dot{x}^{(k+1)},\dot{\lambda}^{(k+1)}),  \label{MSA_lam}\\
	u^{(k+1)} &= \arg\max_{u(t)\in\U} \tilde{H}(x^{(k+1)},\lambda^{(k+1)},u,\dot{x}^{(k+1)},\dot{\lambda}^{(k+1)}). \label{MSA_u}
\end{align}
It is important to note that along solutions to \eqref{PMP_x} and \eqref{PMP_lam}, the right two terms of \eqref{augLag} are zero.  Consequently, only \eqref{MSA_u} is modified with respect to \eqref{PMP_u}.  However, Li et al.~show that this modification is sufficient to ensure convergence \cite{LiChTaE18}.

Li et al.~introduce the regularized forward-backward sweep iteration to train deep neural networks \cite{LiChTaE18} and argue that an advantage of this approach is that it is suitable for application to high dimensional systems.

The analysis of \cite{LiChTaE18} addresses only the continuous time case.  Li et al.~point out that the question of whether Pontryagin's principle holds under numerical discretization is `a delicate one' and refer to counterexamples.  In this paper we show that for variational/symplectic RK methods, an analysis analogous to that of Li et al.~holds.  In particular, their proof of convergence may be translated directly to discrete form.

\section{Variational integrators and symplectic Runge-Kutta pairs\label{sec:VarInt}}

Symplectic Runge-Kutta methods possess two properties that make them attractive for numerical integration of Hamiltonian \emph{initial} value problems: they conserve certain quadratic first integrals and they conserve a modified Hamiltonian function over exponentially long time intervals.  See the monographs \cite{SaSeCa94,HaLuWa06,LeRe05} for a complete discussion.  Symplectic Runge-Kutta methods can be derived using a discrete variational formalism, see \cite{MaWe01}.

Variational methods are also well known in the optimal control literature see e.g.~the work of Marsden, Leok and Ober-Bl\"obaum \cite{ObBlThesis} and references therein. In a recent review, Sanz-Serna \cite{SanzSerna16} argues that it is the property of conservation of quadratic integrals that it is most relevant in the adjoint context.

For optimal control, the use of the variational integrator framework may have additional advantages:  first, by discretizing the integral before optimizing, one constructs a discrete problem for which an optimum may be established, whereas directly discretizing the Euler-Lagrange equations relies on the approximation property in the limit  $\dt\to 0$, where $\dt>0$ is the step size, to guarantee an optimum.  Second, backward error analysis implies the existence of a modified Hamiltonian, near the continuous Hamiltonian, which may have consequences for optimality in the presence of nonunique minima. Backward error analysis may also be applicable for control problems on long time intervals, or for problems with multiple time scales for which the time interval is long on a fast time scale.

We discretize the interval $\T$  into $N>0$ equal steps of size $\dt=T/N$.  An $s$-stage Runge-Kutta method for the state equation \eqref{dynamics} is
\begin{align}
	x_{n+1} &= x_n + \dt \sum_{i=1}^s b_i f(X_{i,n}, U_{i,n}), \label{RK_x} \\
	X_{i,n} &= x_n + \dt \sum_{j=1}^s a_{ij}  f(X_{j,n}, U_{j,n}), \quad i=1,\dots,s, \label{RK_X}
\end{align}
where $n=0,\dots,N-1$ denotes the time step index and the coefficients $b_i$ and $a_{ij}$, $i, j=1,\dots,s$, are chosen to ensure accuracy, stability, and additional properties. See the monographs \cite{HaNoWa,HaWa} for a thorough treatment.  Numerical consistency requires the coefficients $b_i$ satisfy $\sum_i b_i = 1$.  In this paper we will also assume that $b_i\ge 0$, $i=1,\dots,s$.

To simplify notation we will frequently suppress the time step index $n$ in the internal stage variables $X_{i,n}$ and $U_{i,n}$. In all formulas the stage variables are evaluated at time level $n$, so there should be no ambiguity.

A variational integrator for the Lagrangian \eqref{Lagrangian} is a quadrature formula consistent with the above RK method. Enforcing the internal stage relations \eqref{RK_X} requires the introduction of additional Lagrange multipliers. The discrete Lagrangian becomes
\begin{multline}\label{RKalt}
	\Ell[\B{x},\B{\lambda},\B{X},\B{u},\B{G}] = \Phi(x_N) + \lambda_0^T(x_0 - \xi)
	 + \dt \sum_{n=0}^{N-1} \left\{ \phantom{\left( \sum_i^n \frac{X_i}{\dt} - \sum_i^s \right)} \right. \\
	\sum_{i=1}^s b_i h( X_i, U_i) + \lambda_{n+1}^T
	\left(\frac{x_{n+1} - x_n}{\dt} - \sum_{i=1}^s b_i f(X_i,U_i) \right)\\
	 \left.
	- \sum_{i=1}^s b_i G_i \cdot \left(X_i - x_n -  \dt\sum_{j=1}^s a_{ij} f(X_j,U_j)  \right) \right\}.
\end{multline}

Here and henceforth we denote $\B{x} = \{ x_n \, | \, n=0,\dots,N \}$, $\B{X} = \{ X_{i,n} \, | \, i=1,\dots,s; n=0, \dots, N-1\}$, etc.  An exception is the control variable, which only appears at internal stage values.  Consequently we may denote $\B{u} = \{ U_{i,n}\, | \, i=1,\dots,s; n=0, \dots, N-1 \}$ without ambiguity. We also denote $u_n = \{ U_{i,n} \, | \, i=1,\dots,s \}$.

The associated discretization of the cost function \eqref{cost} is
\begin{equation}\label{Jtau}
	J^\dt[\B{u}] = \Phi(x_N) + \dt \sum_{n=0}^{N-1} \sum_{i=1}^s b_i h( X_i, U_i).
\end{equation}
One can formally construct a discrete variational derivative of \eqref{RKalt} with respect to discrete function spaces and a discrete inner product. However for uniform time step $\dt$ it is sufficient to consider just partial derivatives of $\Ell$.  The Euler-Lagrange equations become:
\begin{align}
	\frac{\partial \Ell}{\partial \lambda_n} &= 0 = x_{n+1} - x_n - \dt \sum_{i=1}^s b_i f(X_i,U_i), \quad x_0 = \xi, \label{dEL_x} \\
	\frac{\partial \Ell}{\partial G_i} &= 0 = X_i - x_n - \dt \sum_{j=1}^s a_{ij} f(X_j,U_j),  \label{dEL_X}\\
	\frac{\partial \Ell}{\partial x_n} &= 0 = -\lambda_{n+1} + \lambda_n + \dt \sum_{i=1}^s b_i G_i, \quad \lambda_N = -\prt{\Phi}{x}(x_N),\label{dEL_lam}\\
	\frac{\partial \Ell}{\partial X_k} &= 0 = b_k \prt{h}{x} (X_k,U_k) - b_k \prt{f}{x}(X_k,U_k)^T \lambda_{n+1} -b_k G_k + \dt \sum_{i=1}^s b_i a_{ik} \prt{f}{x} (X_k,U_k)^T G_i,\label{dEL_G}\\
	\frac{\partial \Ell}{\partial U_k} &= 0 = b_k
	\prt{h}{u}(X_k,U_k) - b_k \prt{f}{u}(X_k,U_k)^T\lambda_{n+1} - \dt \sum_{i=1}^s b_i a_{ik} \prt{f}{u}(X_k,U_k)^T G_i. \label{dEL_U}
\end{align}
The relations \eqref{dEL_x}--\eqref{dEL_X} are clearly equivalent to \eqref{RK_x}--\eqref{RK_X}.
Solving \eqref{dEL_lam} for $\lambda_{n+1}$, substituting into \eqref{dEL_G} and defining the coefficients $\tilde{a}_{ij} = b_j - b_j a_{ji}/b_i$,  one finds
\[
	G_i = -\prt{f}{x}(X_i,U_i)^T \left[ \lambda_n + \dt \sum_{j=1}^s \tilde{a}_{ij} G_j \right] + \prt{h}{x} (X_i,U_i).
\]
Similarly \eqref{dEL_U} is written
\begin{equation}\label{Ucond}
	0 = \prt{h}{u}(X_i,U_i) - \prt{f}{u}(X_i,U_i)^T \left[ \lambda_n + \dt \sum_{j=1}^s \tilde{a}_{ij} G_j \right].
\end{equation}
It is useful to introduce the auxiliary stage variable $\Lambda_i$ to represent the term in square brackets in the previous two expressions:
\[
	\Lambda_i = \lambda_n + \dt \sum_{i=1}^s \tilde{a}_{ij} G_j,
\]
such that (cf.~\eqref{gfun})
\[
	G_i = g(X_i,\Lambda_i,U_i) = - \prt{f}{x}(X_i,U_i)^T \Lambda_i + \prt{h}{x}(X_i,U_i)
\]
and the condition \eqref{Ucond} becomes
\[
	0 = \prt{h}{u}(X_i,U_i) - \prt{f}{u}(X_i,U_i)^T\Lambda_i.
\]
In terms of the new variable, the variational Runge-Kutta discretization of Pontryagin's maximum principle is
\begin{align}
	x_{n+1} &= x_n + \dt \sum_{i=1}^s b_i f(X_i, U_i), \quad x_0 = \xi, \label{sum_x}\\
	X_i &= x_n + \dt \sum_{j=1}^s a_{ij}  f(X_j, U_j), \quad i=1,\dots,s, \label{sum_X} \\
	\lambda_{n+1} &= \lambda_n + \dt \sum_{i=1}^s b_i g(X_i,\Lambda_i,U_i), \quad \lambda_N = -\prt{\Phi}{x}(x_N), \label{sum_lam} \\
	 \Lambda_i &= \lambda_{n} + \dt \sum_{j=1}^s \tilde{a}_{ij} g(X_j,\Lambda_j,U_j),\quad i=1,\dots,s,  \label{sum_Lam}\\
	 0 &= \prt{h}{u}(X_i,U_i) - \prt{f}{u}(X_i,U_i)^T\Lambda_i, \quad i=1,\dots,s. \label{sum_U}
\end{align}
This system consists of the state equations \eqref{sum_x} and \eqref{sum_X},
the adjoint equations \eqref{sum_lam} and \eqref{sum_Lam}, and the optimality condition \eqref{sum_U}.

Recalling the Hamiltonian \eqref{Hamiltonian}, we can also write the above relations in a form that emphasizes the Hamiltonian structure:
\begin{align}
	x_{n+1} &= x_n + \dt \sum_{i=1}^s b_i \prt{H}{\lambda}(X_i, \Lambda_i,U_i), \quad x_0 = \xi,\label{RKham_x} \\
	X_i &= x_n + \dt \sum_{j=1}^s a_{ij}  \prt{H}{\lambda}(X_j, \Lambda_j,U_j), \quad i=1,\dots,s, \label{RKham_X} \\
	\lambda_{n+1} &= \lambda_n - \dt \sum_{i=1}^s b_i \prt{H}{x} (X_i,\Lambda_i,U_i), \quad \lambda_N = -\prt{\Phi}{x}(x_N), \label{RKham_lam} \\
	 \Lambda_i &= \lambda_{n} - \dt \sum_{j=1}^s \tilde{a}_{ij} \prt{H}{x} (X_j,\Lambda_j,U_j), \quad i=1,\dots,s, \label{RK_Lam} \\
	 0 &= \prt{H}{u} (X_i,\Lambda_i,U_i), \quad i=1,\dots,s. \label{RKham_U}
\end{align}
In some cases, it is appropriate to replace the latter condition by the more general
\begin{equation}\label{RK_pmp}
	U_i = \arg\max_{u\in\mathcal U} H(X_i,\Lambda_i,u), \quad i=1,\dots,s.
\end{equation}

As noted in \cite{SanzSerna16}, a pair of RK methods defined by coefficients $\{b_i, a_{ij} \}$ and $\{b_i, \tilde{a}_{ij}\}$, where $\tilde {a}_{ij} = b_j - b_j a_{ij}/b_i$, constitute a symplectic partitioned RK pair.  That is, if these methods are applied to a pair of differential equations $\dot x = \prt{H}{\lambda}(x,\lambda)$, $\dot\lambda = -\prt{H}{x}(x,\lambda)$, then the resulting map from $t_n$ to $t_{n+1}$ is a symplectic map.  Hence, we obtain the well-known result that the discrete variational approach automatically produces a symplectic integrator for the Euler-Lagrange equations.

\subsection{Symplectic Euler method}
The elementary example of a symplectic variational integrator is the symplectic Euler method, which
corresponds to the RK pair with $s=1$, $b_1 =1$, $a_{11} = 0 = 1- \tilde{a}_{11}$.  In this case all
 the internal stage relations can be eliminated, leaving the discrete Lagrangian
\begin{equation}
	\Ell[\B{x},\B{\lambda},\B{u}] = \Phi(x_N) + \lambda_0^T(x_0 - \xi) + \dt \sum_{n=0}^{N-1} h(x_n,u_n) + \lambda_{n+1}^T\left( \frac{x_{n+1}-x_n}{\dt} - f(x_n,u_n) \right).
\end{equation}
The discrete Pontryagin maximum principle is
\begin{align}
	x_{n+1} &= x_n + \dt f(x_n,u_n),  \label{Euler_x} \\
	\lambda_{n+1} &= \lambda_n - \dt \prt{f}{x}(x_n,u_n)^T\lambda_{n+1} + \dt \prt{h}{x}(x_n,u_n), \label{Euler_lam}  \\
	0 &= \prt{f}{u}(x_n,u_n)^T\lambda_{n+1} - \prt{h}{u}(x_n,u_n), \label{Euler_u}
\end{align}
with boundary conditions $x_0 = \xi$, $\lambda_N = -\prt{\Phi}{x}(x_N)$.

Note that \eqref{Euler_x}--\eqref{Euler_u} can also be written in terms of the Hamiltonian $H$:
\begin{align}
	\frac{x_{n+1}-x_n}{\dt} &= \prt{H}{\lambda}(x_n, \lambda_{n+1},u_n), \label{Hager_SE1}\\
	\frac{\lambda_{n+1} - \lambda_n}{\dt} &= -\prt{H}{x}(x_n, \lambda_{n+1},u_n), \label{Hager_SE2}\\
	0 &= \prt{H}{u}(x_n,\lambda_{n+1},u_n). \label{Hager_SE3}
\end{align}

\subsection{Reduced notation for Runge-Kutta methods}
Hager \cite{Hager00} introduced notation that casts general symplectic Runge-Kutta methods \eqref{RKham_x}--\eqref{RKham_U} in a form consistent with the symplectic Euler method.
Define
\begin{equation}\label{ftau}
	f^\dt(x,u) = \sum_{i=1}^s b_i f(X_i(x,u),U_i(u)), \quad
	h^\dt(x,u) = \sum_{i=1}^s b_i h(X_i(x,u),U_i(u)),
\end{equation}
where we view the stage values $X_i$ and $U_i$ as functions of grid point value $x$ and discrete control  $u = \{ U_1,\dots, U_s \}$ according to
\begin{equation}\label{Xmap}
	X_i(x,u) = x + \dt \sum_{j=1}^s a_{ij} f(X_j(x,u),U_j(u)), \quad i =1,\dots,s.
\end{equation}
Similarly, define the Hamiltonian
\begin{equation}\label{Htau}
	H^\dt(x,\lambda,u) = \lambda^T f^\dt(x,u) - h^\dt(x,u).
\end{equation}
With this notation, the discretization of Pontryagin's maximum principle with any symplectic Runge-Kutta pair can be written as
\begin{align}
	\frac{x_{n+1} - x_n}{\dt} &= \prt{H^\dt}{\lambda}(x_n,\lambda_{n+1},u_n), \label{Hager_x}\\
	\frac{\lambda_{n+1}-\lambda_n}{\dt} &= -\prt{H^\dt}{x} (x_n,\lambda_{n+1},u_n),\label{Hager_lam}\\
	0 &= \prt{H^\dt}{u}(x_n,\lambda_{n+1},u_n).\label{Hager_u}
\end{align}
To see the equivalence, note that evaluating \eqref{Xmap} at $x_n$ yields the implicit relations \eqref{sum_X}. Taking the derivative of \eqref{Htau} with respect to $\lambda$ and substituting \eqref{ftau} shows \eqref{Hager_x} to be equivalent to \eqref{sum_x}.  The proof of the relation \eqref{Hager_lam} is more involved. We adapt the proof from \cite{Hager00} to our notation.

Let $\Psi_i(x) = \partial_x X_i(x,u)$ and denote $\Psi_i = \Psi_i(x_n)$. Then computing the derivative of \eqref{Xmap} at $x_n$ yields the linear system
\begin{equation}\label{Psi}
	\Psi_i = I + \dt \sum_j a_{ij} \prt{f}{x}(X_i,U_i)  \Psi_j.
\end{equation}
The derivative on the right side of \eqref{Hager_lam} is
\begin{equation}\label{Hager_Hx}
	\prt{H^\dt}{x} (x_n,\lambda_{n+1},u_n) = \sum_{j=1}^s  b_j \Psi_j^T\prt{f}{x}(X_j,U_j)^T \lambda_{n+1}  - b_j  \Psi_j^T \prt{h}{x} (X_j,U_j).
\end{equation}
Rearranging \eqref{dEL_G} gives
\[
 		b_j G_j - \dt \sum_{i=1}^s b_i a_{ij} \prt{f}{x} (X_j,U_j)^T G_i = b_j \prt{h}{x} (X_j,U_j) - b_j \prt{f}{x}(X_j,U_j)^T \lambda_{n+1}.
\]
Premultiplying by $\Psi_j^T$ and summing over $j$ gives
\begin{multline}
	\sum_{j=1}^s b_j \Psi_j^TG_j - \dt \sum_{i,j=1}^s b_i a_{ij} \Psi_j^T\prt{f}{x} (X_j,U_j)^TG_i\\
	= \sum_{j=1}^s b_j \Psi_j^T\prt{h}{x} (X_j,U_j)  - b_j  \Psi_j ^T\prt{f}{x}(X_j,U_j)^T \lambda_{n+1}= -\prt{H^\dt}{x} (x_n,\lambda_{n+1},u_n),
\end{multline}
where the last equality follows from \eqref{Hager_Hx}.
Now changing the index of summation in the first sum on the left, we obtain
\begin{align*}
	-\prt{H^\dt}{x} (x_n,\lambda_{n+1},u_n) &= \sum_{i=1}^s b_i \Psi_i^T G_i  - \dt \sum_{i=1}^s \left( \sum_{j=1}^s
	a_{ij} \Psi_j^T \prt{f}{x} (X_j,U_j)^T\right) b_i G_i\\
	&= \sum_{i=1}^s b_i G_i\\
	&= \frac{\lambda_{n+1} - \lambda_n}{\dt},
\end{align*}
where the second equality follows from \eqref{Psi}, thus confirming \eqref{Hager_lam}.

The proof of \eqref{RKham_U} follows similar arguments, see \cite{Hager00}.
Note the analogy between the relations \eqref{RKham_x}--\eqref{RKham_U} and \eqref{Hager_SE1}--\eqref{Hager_SE2} for the symplectic Euler method.

\section{Convergence analysis\label{sec:convergence}}
In this section we prove the convergence of the regularized forward-backward sweep iteration \eqref{MSA_x}--\eqref{MSA_u} for symplectic Runge-Kutta methods. The proof here follows closely that of Li et al.~for the continuous case \cite{LiChTaE18}.  It is the symplectic/variational structure that facilitates this analogy.

Using the compact notation \eqref{ftau} and \eqref{Htau}, we define the discrete regularized Hamiltonian function
\begin{equation}\label{Htilde_tau}
	\tilde{H}^\dt(x,\lambda,u,q,p) = H^\dt(x,\lambda,u) - \frac{\rho}{2}\left( \| q - \prt{H^\dt}{\lambda}(x,\lambda,u) \|^2 + \| p + \prt{H^\dt}{x} (x,\lambda,u) \|^2\right).
\end{equation}
In iterate $k$, the symplectic Runge-Kutta discretization of the regularized forward-backward sweep iteration \eqref{MSA_x}--\eqref{MSA_u}   solves, in sequence,
\begin{align}
	x_{n+1}^{(k+1)} &= x_n^{(k+1)}  + \dt \prt{\tilde{H}^\dt}{\lambda} \left(x_n^{(k+1)} ,\lambda_{n+1}^{(k)},u_n^{(k)},\frac{x_{n+1}^{(k+1)} -x_n^{(k+1)} }{\dt},\frac{\lambda_{n+1}^{(k)}-\lambda_n^{(k)}}{\dt}\right),  \label{dMSA_x} \\
	\lambda_{n+1}^{(k+1)}  &= \lambda_n^{(k+1)}  - \dt\prt{\tilde{H}^\dt}{x} \left(x_n^{(k+1)} ,\lambda_{n+1}^{(k+1)} ,u_n^{(k)},\frac{x_{n+1}^{(k+1)} -x_n^{(k+1)} }{\dt},\frac{\lambda_{n+1}^{(k+1)} -\lambda_n^{(k+1)} }{\dt}\right), \label{dMSA_lam}  \\
	u_n^{(k+1)}  &= \arg\max_{u\in\U} \tilde{H}^\dt\left(x_n^{(k+1)} ,\lambda_{n+1}^{(k+1)} ,u,\frac{x_{n+1}^{(k+1)} -x_n^{(k+1)} }{\dt},\frac{\lambda_{n+1}^{(k+1)} -\lambda_n^{(k+1)} }{\dt}\right),  \label{dMSA_u}
\end{align}
proceeding as follows: \eqref{dMSA_x} by forward integration with $\B{u}$ and $\B{\lambda}$ fixed, then \eqref{dMSA_lam} by backward integration with $\B{x}$ and $\B{u}$ fixed, and finally \eqref{dMSA_u} solved for each time step independently (e.g.~in parallel), with $\B{x}$ and $\B{\lambda}$ fixed.

It is important to recall that with $u$ fixed, along solutions of \eqref{dMSA_x} and \eqref{dMSA_lam} the extra regularization terms in the extended Hamiltonian $\tilde{H}^\dt$ are identically zero and
\begin{align*}
	\prt{\tilde{H}^\dt}{\lambda}  \left(x_n,\lambda_{n+1},u_n,\frac{x_{n+1}-x_n}{\dt},\frac{\lambda_{n+1}-\lambda_n}{\dt}\right) &= \prt{H^\dt}{\lambda}(x_n,\lambda_{n+1},u_n), \\
	\prt{\tilde{H}^\dt}{\lambda} \left(x_n,\lambda_{n+1},u_n,\frac{x_{n+1}-x_n}{\dt},\frac{\lambda_{n+1}-\lambda_n}{\dt}\right) &= \prt{H^\dt}{x}(x_n,\lambda_{n+1},u_n),
\end{align*}
 i.e., the regularization terms only affect the maximization step \eqref{dMSA_u}.

\subsubsection*{Notation and identities}

In the following we consider a single iteration of \eqref{dMSA_x}--\eqref{dMSA_u}. We think of $H^\dt$, $\B{x}$ and $\B{\lambda}$ as functions of $\B{u}$.  Consequently we denote by $x_n^u$  and $\lambda_n^u$ the numerical solutions to \eqref{Hager_x} and \eqref{Hager_lam} given a candidate control $\B{u}$.

It is convenient to define the composite notation
\[
	 z_n = \begin{pmatrix} x_n \\ \lambda_{n+1} \end{pmatrix}, \quad \prt{H^\dt}{z}(z_n,u_n) = \begin{pmatrix} \prt{H^\dt}{x}(x_n,\lambda_{n+1},u_n) \\ \prt{H^\dt}{\lambda}(x_n,\lambda_{n+1},u_n) \end{pmatrix}.
\]

We consider two control sequences $\B{u}$ and $\B{v}$, and we are interested in bounding the change in $\tilde{H}^\tau$ when $\B{u}$ is replaced by $\B{v}$. To that end we define an operator that denotes the difference between quantities dependent on $\B{u}$ and $\B{v}$:
\[
	\Du x_n = x_n^v - x_n^u.
\]
We use this notation also for functions, e.g.
\[
	\Du H^\dt|_n = H^\dt(z_n^v,v_n) - H^\dt(z_n^u,u_n).
\]
We denote by $\Dbar H^\dt$ the change due to an update in $\B{u}$ with $\B{x}$ and $\B{\lambda}$ fixed as functions of $u$:
\begin{equation}\label{Ham_incr}
	\Dbar H^\dt|_n = H^\dt(x_n^u,\lambda_{n+1}^u,v_n) - H^\dt(x_n^u,\lambda_{n+1}^u,u_n).
\end{equation}
We denote the temporal forward difference operator by $\FD$:
\[
	\FD x_n = \frac{x_{n+1} - x_n}{\dt},
\]
and remark that $\FD$ commutes with $\Du$ when applied to variables, i.e. $\Du \FD x_n =  \FD\Du x_n$.

Next we note the discrete integration by parts formula:
\begin{align*}
	\dt \sum_{n=0}^{N-1} \lambda_{n+1}^T \FD x_n &=  \sum_{n=0}^{N-1} \lambda_{n+1}^T (x_{n+1} - x_n) \\
	&= -\lambda_0^T x_0 + \lambda_0^T x_0 - \lambda_1^T x_0 + \lambda_1^T x_1 + \cdots - \lambda_N^T x_{N-1} + \lambda_N^T x_N \\
	&= \lambda_n^T x_n\big|_0^N - \dt \sum_{n=0}^{N-1} (\FD \lambda_n)^T x_n.
\end{align*}
This formula holds for any discrete functions defined for $n=0,\dots,N$, and in particular we may insert the difference operator $\Du$ to obtain two useful alternatives:
\begin{equation} \label{parts1}
	\dt \sum_{n=0}^{N-1} \lambda_{n+1}^u \cdot \FD \Du x_n = \lambda_n^u \cdot \Du x_n \big|_0^N
	- \dt \sum_{n=0}^{N-1} \FD \lambda_n^u\cdot \Du x_n,
\end{equation}
\begin{equation} \label{parts2}
	\dt \sum_{n=0}^{N-1} \Du\lambda_{n+1} \cdot \FD \Du x_n = \Du \lambda_n \cdot \Du x_n \big|_0^N
	- \dt \sum_{n=0}^{N-1} \FD \Du \lambda_n \cdot \Du x_n.
\end{equation}

\subsubsection*{Estimates}
In the Appendix we show that---possibly with a restriction on step size---the Lipschitz conditions \eqref{contLipschitz} on $f$ and $h$ translate into related Lipschitz conditions on $f^\dt$ and $h^\dt$.  Henceforth choosing $K$ to be a generic Lipschitz constant we obtain the bounds
\begin{equation} \label{discLipschitz} \begin{split}
	&\| f^\dt(x,u) - f^\dt(x',u) \| + \| \prt{f^\dt}{x}(x,u) -  \prt{f^\dt}{x}(x',u) \| \le K \| x - x'\|, \\
	&| h^\dt(x,u) - h^\dt(x',u) | + \| \prt{h^\dt}{x}(x,u) -  \prt{h^\dt}{x}(x',u) \| \le K \| x - x'\|.
	\end{split}
\end{equation}
Note also that the leftmost terms in the above inequalities as well as the analogous ones of \eqref{contLipschitz} imply global bounds on the derivatives (which may be relaxed, see \cite{LiChTaE18})
\begin{equation} \label{bound}
	\| \prt{\Phi}{x}(x) \| \le K, \quad
	\| \prt{f}{x}(x,u) \| \le K, \quad
	\| \prt{h}{x}(x,u) \| \le K, \quad
	\| \prt{f^\dt}{x}(x,u) \| \le K, \quad
	\| \prt{h^\dt}{x}(x,u) \| \le K.
\end{equation}

We use two discrete forms of Gr\"onwall's lemma \cite{Emmrich99}.  Let $\{b_n\}$ be a given, monotone sequence and $\dt,K>0$.  Then the following implication holds:
\begin{equation} \label{gronwall1}
	a_{n+1} \le (1+\dt K) a_n + \dt b_n, \quad \forall n \quad \Rightarrow \quad
	a_{n} \le e^{\dt n K} a_0 + K^{-1} e^{\dt n K} b_{n-1}.
\end{equation}
Under the same conditions, the following implication holds:
\begin{equation}\label{gronwall2}
	a_{n+1} \le b_{n+1} + \dt K \sum_{m=0}^n a_m, \quad \forall n \quad \Rightarrow \quad
	a_n \le e^{\dt n K } b_n.
\end{equation}

From \eqref{Hager_Hx} and \eqref{Psibound}, and using the bounds \eqref{bound} on $\prt{f}{x}$ and $\prt{h}{x}$,
\[
	\| \lambda_n \| \le \| \lambda_{n+1} \|  + \dt \| \prt{H^\dt}{x} (x_n \lambda_{n+1}, u_n) \|
	\le (1 + \dt K ) \| \lambda_{n+1} \| + \dt K,
\]
where we have absorbed the constant from \eqref{Psibound} into $K$.
Further using Gr\"onwall bound \eqref{gronwall1} and the bound \eqref{bound} on $\prt{\Phi}{x}(x)$,
\begin{equation} \label{lam_bound}
	\| \lambda_n \| \le K_1 := (K+1) e^{\dt K N} = (K+1)e^{KT}.
\end{equation}

From $\Du x_{n+1} = \Du x_n + \dt \Du f^\dt|_n$  and $\Du x_0=0$ we calculate
\begin{align*}
	\| \Du x_n \| &\le \dt \sum_{m=0}^{n-1} \| \Du f^\dt|_m \| \\
	&\le \dt \sum_{m=0}^{n-1}\| \Dbar f^\dt|_m \| + \| f^\dt(x_m^v,v_m) - f^\dt(x_m^u,v_m) \| \\
	&\le \dt \sum_{m=0}^{n-1}\| \Dbar f^\dt|_m \| + K \| \Du x_m \|,
\end{align*}
and using Gr\"onwall bound \eqref{gronwall2},
\begin{equation}\label{Du_x_bound}
	\| \Du x_n \| \le \dt e^{K T} \sum_{m=0}^{N-1} \| \Dbar f^\dt|_m \|.
\end{equation}

Similarly, from $\Du \lambda_n = \Du \lambda_{n+1} + \dt \Du \prt{H^\dt}{x}(x_n,\lambda_{n+1},u_n)$ we obtain
\begin{align*}
	\| \Du \lambda_n \| &\le \| \Du \lambda_N \| + \dt \sum_{m=n}^{N-1} \| \Du \prt{H^\dt}{x}|_m \| \\
	&\le K \| \Du x_N \| + \dt \sum_{m=n}^{N-1} \| \Dbar \prt{H^\dt}{x}|_m \|  + \dt K \sum_{m=n}^{N-1} \| \Du \lambda_{m+1} \|
		+ \dt K (K_1+1) \sum_{m=n}^{N-1} \| \Du x_m\|,
\end{align*}
where the last term uses \eqref{contLipschitz} and the Lipschitz condition \eqref{discLipschitz} on $\prt{H^\dt}{x}$.  The discrete
Gr\"{o}nwall's lemma gives
\[
	\| \Du \lambda_n \| \le K e^{K T} \left( \| \Du x_N \| + \dt (K_1+1)  \sum_{m=0}^{N-1} \| \Du x_m \| \right)  + \dt e^{K T} \sum_{m=0}^{N-1} \| \Dbar \prt{H^\dt}{x}|_m \|.
\]
Finally, making use of \eqref{Du_x_bound} gives
\begin{equation} \label{Du_lam_bound}
\| \Du \lambda_{n} \| \le \dt  K_2 \sum_{m=0}^{N-1} \| \Dbar f^\dt|_m \| + \dt e^{K T}  \sum_{m=0}^{N-1} \| \Dbar \prt{H^\dt}{x}|_m \|, \quad K_2 = K e^{2K T} (1 + (K_1+1)T).
\end{equation}

The following estimates make use of Taylor's theorem in the mean value form:
\begin{equation}\label{DuHz}
	\Du \prt{H^\dt}{z}|_n \cdot \Du z_n = \Dbar \prt{H^\dt}{z}|_n \cdot \Du z_n + \Du z_n\cdot \prt{H^\dt}{zz} ( z_n^u + r_1 \Du z_n,u_n) \cdot \Du z_n,
\end{equation}
for some $r_1 \in [0,1]$, where $\prt{H^\dt}{zz}$ denotes the Hessian matrix of second partial derivatives of $H^\dt$.
\begin{equation} \label{DuPsix}
	\Du \prt{\Phi}{x}(x_N) \cdot \Du x_N = \Du x_N\cdot \prt{\Phi}{xx} (x_N^u + r_2 \Du x_N)\cdot \Du x_N,
\end{equation}
for some $r_2\in [0,1]$. Similarly,
\begin{equation}\label{DuPsi}
	\prt{\Phi}{x}(x_N^u)\cdot \Du x_N = \Phi(x_N^v) - \Phi(x_N^u) - \frac{1}{2} \Du x_N\cdot \prt{\Phi}{xx} (x_N^u + r_3 \Du x_N)\cdot \Du x_N,
\end{equation}
for some $r_3\in [0,1]$.
\begin{equation}\label{DuH}
	\Du H^\dt = \Dbar H^\dt + \prt{H^\dt}{z}(z_n^u,v)\cdot \Du z_n  + \frac{1}{2} \Du z_n \cdot \prt{H^\dt}{zz}(z_n^u+r_4\Du z_n,v_n)\cdot \Du z_n,
\end{equation}
for some $r_4\in [0,1]$.

\subsubsection*{Convergence of the iteration}

Convergence of the regularized forward-backward sweep iteration relies on Lemma 2 of \cite{LiChTaE18}, the proof of which we adapt for the symplectic RK method here.  The result we want states that under the assumptions \eqref{contLipschitz}, there exists a constant $C>0$ such that for any two discrete controls $\B{u},\B{v} \in \U$, the discrete cost function \eqref{Jtau} satisfies
\begin{multline} \label{Lemma2}
	J^\dt(\B{v}) \le J^\dt(\B{u}) - \dt \sum_{n=0}^{N-1} \Dbar H^\dt|_n + C \dt \sum_{n=0}^{N-1} \| f^\dt(x_n^u,v_n) - f^\dt(x_n^u,u_n) \|^2  \\
	+ C \dt \sum_{n=0}^{N-1} \| \prt{H^\dt}{x}(x_n^u,\lambda_{n+1}^u,v_n) - \prt{H^\dt}{x}(x_n^u,\lambda_{n+1}^u,u_n) \|^2 \\
	= J^\dt(\B{u}) - \dt \sum_{n=0}^{N-1} \Dbar H^\dt|_n + C\dt \sum_{n=0}^{N-1} \| \Dbar \prt{H^\dt}{z}|_n \|^2.
\end{multline}

Define the discrete functional
\begin{equation} \label{Idef}
	\mathcal{I} (\B{x},\B{\lambda},\B{u}) = \dt \sum_{n=0}^{N-1} \lambda_{n+1}^T \FD x_n - H^\dt(x_n,\lambda_{n+1},u_n) - h^\dt(x_n,u_n) \equiv 0.
\end{equation}
The functional $\mathcal{I}$ is identically zero for sequences $\B{x}$ and $\B{\lambda}$ satisfying \eqref{Hager_x}--\eqref{Hager_lam}.
Note the identity
\begin{equation}\label{id-1form}
	\Du (\lambda_{n+1}\cdot \FD x_n) =
	\lambda_{n+1}^{u} \cdot \FD \Du x_n + \Du \lambda_{n+1} \cdot \FD x_n^u + \Du \lambda_{n+1} \cdot \FD \Du x_n,
\end{equation}

we find
\begin{multline*}
	0 \equiv \mathcal{I}(\B{x}^v,\B{\lambda}^v,\B{v}) - \mathcal{I}(\B{x}^u,\B{\lambda}^u, \B{u}) = \\
	\dt \sum_{n=0}^{N-1} \lambda_{n+1}^{u} \cdot \FD\Du x_n + \Du \lambda_{n+1} \cdot \FD x_n^u +
	\Du \lambda_{n+1} \cdot \FD \Du x_n \\
	- \dt \sum_{n=0}^{N-1} \left(H^\dt(x_n^v,\lambda_{n+1}^v,v_n) - H^\dt(x_n^u,\lambda_{n+1}^u,u_n)\right) \\
	- \dt \sum_{n=0}^{N-1} \left( h^\dt(x_n^v,v_n) - h^\dt(x_n^u,u_n) \right).
\end{multline*}
In our notation this is
\begin{equation}
	0 \equiv \Du \mathcal{I} = \dt \sum_{n=0}^{N-1} \lambda_{n+1}^{u} \cdot \FD\Du x_n + \Du \lambda_{n+1} \cdot \FD x_n^u +
	\Du \lambda_{n+1} \cdot \FD \Du x_n
	-  \Du H^\dt|_n - \Du h^\dt|_n. \label{dI}
\end{equation}
{\bf Remark.} {\it  This is the point where the symplectic/variational property of the symplectic RK method is important.
Since $x_n$ and $\lambda_n$ are discretized by a symplectic partitioned Runge-Kutta method, we see that $\mathcal{I}$ is also equivalent to
the constraint part of the discrete Lagrangian:
\[
	\mathcal{I} = \dt\sum_{N=0}^{N-1} \lambda_{n+1}^T\left( \frac{x_{n+1}-x_n}{\dt} - f^\dt(x_n,u_n)\right),
\]
which is identically zero along a solution to the state dynamics \eqref{Hager_x}.  Of course, one could define $\mathcal{I}$ as above for an arbitrary choice of the $\lambda_n$.  Then $\mathcal{I}$ would be identically zero, but one would not be able to translate this into a statement about the Hamiltonian.}

Using \eqref{parts1} the first two terms on the right side of \eqref{dI} are equal to
\begin{multline*}
	\dt \sum_{n=0}^{N-1}  \lambda_{n+1}^{u} \cdot \FD\Du x_n + \Du \lambda_{n+1} \cdot \FD x_n^u \\
	= \lambda_n^u \cdot \Du x_n \big|_0^N
	+ \dt \sum_{n=0}^{N-1}  f^\dt(x_n^u,u_n) \cdot \Du \lambda_{n+1} + \prt{H^\dt}{x} (x_n^u,\lambda_{n+1}^u, u_n)\cdot \Du x_n,
\end{multline*}
or in compact notation
\begin{equation}\label{dI_terms12}
	\dt \sum_{n=0}^{N-1}  \lambda_{n+1}^{u} \cdot \FD\Du x_n + \Du \lambda_{n+1} \cdot \FD x_n^u
	= \lambda_n^u \cdot \Du x_n \big|_0^N + \dt \sum_{n=0}^{N-1} \prt{H^\dt}{z}(z_n^u,u_n) \cdot \Du z_n.
\end{equation}
Similarly, using \eqref{parts2} the third term on the right side of \eqref{dI} is equal to
\begin{multline*}
	\dt \sum_{n=0}^{N-1} \Du \lambda_{n+1} \cdot  \FD \Du x_n
	= \frac{1}{2} \dt \sum_{n=0}^{N-1} \Du \lambda_{n+1} \cdot \FD\Du x_n +
	\frac{1}{2} \dt \sum_{n=0}^{N-1} \Du \lambda_{n+1} \cdot \FD\Du x_n \\
	= \frac{1}{2}\Du \lambda_n \cdot \Du x_n\big|_0^N + \frac{1}{2}
	\dt\sum_{n=0}^{N-1} \left( \prt{H^\dt}{x}(x_n^v,\lambda_{n+1}^v,v_n) - \prt{H^\dt}{x} (x_n^u,\lambda_{n+1}^u,u_n)\right) \cdot \Du x_n \\
	+
	\left( \prt{H^\dt}{\lambda}(x_n^v,\lambda_{n+1}^v,v_n) - \prt{H^\dt}{\lambda} (x_n^u,\lambda_{n+1}^u,u_n)\right) \cdot \Du \lambda_{n+1},
\end{multline*}
or,
\begin{equation}\label{dI_term3}
	\dt \sum_{n=0}^{N-1} \Du \lambda_{n+1} \cdot  \FD \Du x_n 	
	= \frac{1}{2}\Du \lambda_n \cdot \Du x_n\big|_0^N
	+ \frac{1}{2} \dt\sum_{n=0}^{N-1} \Du \prt{H^\dt}{z}|_n \cdot \Du z_n.
\end{equation}
{\bf Remark.}  {\it Again the symplectic property of the discretization allows us to express this as the gradient of the Hamiltonian collocated at the numerical solution of the forward and backward equations, which in turn will allow cancellation with the second term of the Taylor expansion in \eqref{dI_simplified}.}

Combining \eqref{dI}, \eqref{dI_terms12} and \eqref{dI_term3} gives
\begin{multline} \label{dIprime}
	0 \equiv \Du \mathcal{I} = (\lambda_n^u  + \frac{1}{2}\Du \lambda_n) \cdot \Du x_n\big|_0^N + \\
	\dt \sum_{n=0}^{N-1}  \prt{H^\dt}{z}(z_n^u,u_n) \cdot \Du z_n  + \frac{1}{2}  \Du \prt{H^\dt}{z}|_n \cdot \Du z_n -  \Du H^\dt|_n - \Du h^\dt|_n.
\end{multline}

Given that $\Du x_0= 0$, the boundary term in \eqref{dIprime} reduces to
\begin{equation}\label{boundary}
	 (\lambda_N^u + \frac{1}{2} \Du \lambda_N) \cdot \Du x_N  = -\prt{\Phi}{x}(x_N)\cdot\Du x_N - \frac{1}{2} \left(\prt{\Phi}{x}(x_N^v) - \prt{\Phi}{x}(x_N^u)\right)\cdot \Du x_N.
\end{equation}

We substitute \eqref{DuHz} and \eqref{DuH} into the second and third summand of \eqref{dIprime}, \eqref{boundary} into the boundary term, and subsequently the estimates \eqref{DuPsix} and \eqref{DuPsi} to yield:
\begin{multline*}
	0 \equiv \Du\mathcal{I}  = -\left( \Phi(x_N^v) - \Phi(x_N^u) - \frac{1}{2} \Du x_N\cdot \prt{\Phi}{xx} (x_N^u + r_3 \Du x_N)\cdot \Du x_N \right) \\
	- \frac{1}{2} \left( \Du x_N\cdot \prt{\Phi}{xx} (x_N^u + r_2 \Du x_N)\cdot \Du x_N\right)
	+ \dt \sum_{n=0}^{N-1}  -\Du h^\dt|_n + \prt{H^\dt}{z}(z_n^u,u_n) \cdot \Du z_n   \\
	+ \frac{1}{2} \left( \Dbar \prt{H^\dt}{z}|_n \cdot \Du z_n + \Du z_n\cdot \prt{H^\dt}{zz} ( z_n^u + r_1 \Du z_n,u_n) \cdot \Du z_n \right) \\
	-\left(\Dbar H^\dt|_n + \prt{H^\dt}{z}(z_n^u,v_n)\cdot \Du z_n  + \frac{1}{2} \Du z_n \cdot \prt{H^\dt}{zz}(z_n^u+r_4\Du z_n,v_n)\cdot \Du z_n\right),
\end{multline*}
or,
\begin{multline}\label{dI_simplified}
\Du \Phi(x_N) + \dt\sum_{n=0}^{N-1} \Du h^\dt(x_n,u_n) = \\
 - \frac{1}{2} \Du x_N\cdot \left( \prt{\Phi}{xx} (x_N^u + r_2 \Du x_N) - \prt{\Phi}{xx} (x_N^u + r_3 \Du x_N)\right) \cdot \Du x_N  \\
 - \dt \sum_{n=0}^{N-1} \Dbar H^\dt|_n + \frac{1}{2} \dt \sum_{n=0}^{N-1} \Dbar \prt{H^\dt}{z}|_n \cdot \Du z_n \\
 + \frac{1}{2}\dt \sum_{n=0}^{N-1}  \Du z_n \cdot \left( \prt{H^\dt}{zz}(z_n^u+r_1\Du z_n,v_n) - \prt{H^\dt}{zz}(z_n^u+r_4\Du z_n,v_n)\right) \cdot \Du z_n.
\end{multline}

Next, we use the estimates \eqref{Du_x_bound} and \eqref{Du_lam_bound} and the fact that the quadratic terms are bounded by some constant $K_3$ to calculate
\begin{align*}
	J^\dt[\B{v}] - J^\dt[\B{u}] \le &-\dt \sum_{n=0}^{N-1} \Dbar H^\dt|_n \\
	&+ K_3 \| \Du x_N\|^2 + K_3 \dt \sum_{n=0}^{N-1} \left(  \| \Du x_n \|^2 + \| \Du \lambda_{n+1} \|^2 \right) \\
	&+ \frac{1}{2} \dt \sum_{n=0}^{N-1} \| \Du x_n\| \| \Dbar f^\dt|_n\| + \frac{1}{2} \dt \sum_{n=0}^{N-1}  \| \Du \lambda_{n+1} \| \| \Dbar \prt{H^\dt}{x}|_n \| \\
    \le & -\dt \sum_{n=0}^{N-1} \Dbar H^\dt|_n + C \left( \dt \sum_{n=0}^{N-1} \| \Dbar f^\dt|_n \| \right)^2 + C \left(\dt \sum_{n=0}^{N-1} \| \Dbar \prt{H^\dt}{x}|_n  \| \right)^2 \\
    \le & -\dt \sum_{n=0}^{N-1} \Dbar H^\dt|_n + C \dt \sum_{n=0}^{N-1} \| \Dbar f^\dt|_n \|^2 + C \dt \sum_{n=0}^{N-1} \| \Dbar  \prt{H^\dt}{x}|_n \|^2,
\end{align*}	
which is the result sought (cf.~\eqref{Lemma2}).
	
It now remains to show that the regularized forward-backward sweep iteration converges.  We first show that an estimate of the same form as \eqref{Lemma2} holds for  $\Du H^\dt$ when the regularized Hamiltonian is maximized.  These can be combined to show monotone decay of the objective function $J^\dt[\B{u}]$.  Thereafter, it is shown that the sum of the decrements is finite, which implies convergence of the differences.


Let $\B{v}$ denote the improved control obtained by solving \eqref{dMSA_u}.  The resulting change in $\tilde{H}^\dt$ must be nonnegative, hence
\begin{multline}
	0 \le \dt  \sum_{n=0}^{N-1}\Dbar \tilde{H}^\dt|_n= \dt  \sum_{n=0}^{N-1} \Dbar H^\dt|_n \\
	- \frac{\rho}{2} \left[  \| \frac{x_{n+1}^u-x_n^u}{\dt} - f^\dt(x_n^u,v_n) \|^2 +
	\| \frac{\lambda_{n+1}^u-\lambda_n^u}{\dt} +H^\dt_x (x_n^u,\lambda_{n+1}^uv_n) \|^2 \right] \\
	 + \frac{\rho}{2} \left[ \| \frac{x_{n+1}^u-x_n^u}{\dt} - f^\dt(x_n^u,u_n) \|^2 +
	\| \frac{\lambda_{n+1}^u-\lambda_n^u}{\dt} +H^\dt_x (x_n^u,\lambda_{n+1}^u,u_n) \|^2 \right].
\end{multline}
The last term in square brackets vanishes since $x_n^u$ and $\lambda_n^u$ satisfy \eqref{Hager_x}--\eqref{Hager_lam}.
Consequently, the above expression is equivalent to
\begin{equation} \label{duHt}
	0 \le \dt  \sum_{n=0}^{N-1}\Dbar \tilde{H}^\dt|_n = \dt  \sum_{n=0}^{N-1} \Dbar H^\dt|_n - \frac{\rho}{2} \left[  \| \Dbar f^\dt|_n \|^2 +
	\| \Dbar \prt{H^\dt}{x}|_n \|^2 \right].
\end{equation}
Combining this with Lemma 2 gives
\begin{equation}\label{monotone}
	J^\dt[\B{v}] - J^\dt[\B{u}] \le -(1- \frac{2C}{\rho}) \dt \sum_{n=0}^{N-1} \Dbar H^\dt|_n.
\end{equation}
The summation on the right side is nonnegative, as a consequence of \eqref{duHt} .  Therefore, choosing $\rho>2C$ ensures that $J^\dt$ is nonincreasing.  Next suppose we iterate \eqref{dMSA_x}--\eqref{dMSA_u}.  Let $\B{u}^{(k)}$ denote the control variable in iteration $k$.  Then it holds that
\[
	\sum_{k=0}^M \dt \sum_{n=0}^{N-1} \Dbar H^\dt|_n^{(k)} \le D^{-1} (J^\dt[\B{u}^{(0)}] - J^\dt[\B{u}^{(M+1)}]) \le D^{-1}( J^\dt[\B{u}^{(0)}] - \inf_{\B{u}\in \U} J^\dt[\B{u}] ),
\]
where $D = (1-2C/\rho)>0$.
Consequently, in the limit $M\to\infty$ this sum is bounded, which implies
\[
	\sum_{n=0}^{N-1} \Dbar H^\dt|_n \to 0,
\]
proving convergence of the iteration.

\section{Numerical illustration\label{sec:numer}}
In this section we study numerically the convergence of the discrete regularized forward-backward sweep iteration. As a test problem we control the motion of a damped oscillator in a double well potential. The controlled motion is given by
\begin{equation}
	x = \begin{pmatrix} q \\ p \end{pmatrix}, \qquad f(x,u) = \begin{pmatrix} p \\ q - q^3 - \nu p + u \end{pmatrix}, \label{dwp}
\end{equation}
where $\nu>0$ is a damping parameter.  The control $u(t)$ acts only on the velocity. As initial condition we choose $\xi = (-1,0)$ in the left potential well, and we seek to minimize the cost function
\begin{equation}
	J[u] = \frac{\alpha}{2} \| x(T) - x_f \|^2 + \int_0^T \frac{1}{2} u(t)^2 \, dt,
\end{equation}
where the target final position is $x_f = (1,0)$, in the right potential well. For the numerical computations we take $T=6$, $\nu = 1$, and $\alpha = 10$.

\begin{figure}[!htbp]
\begin{center}
\includegraphics[width=0.49\textwidth]{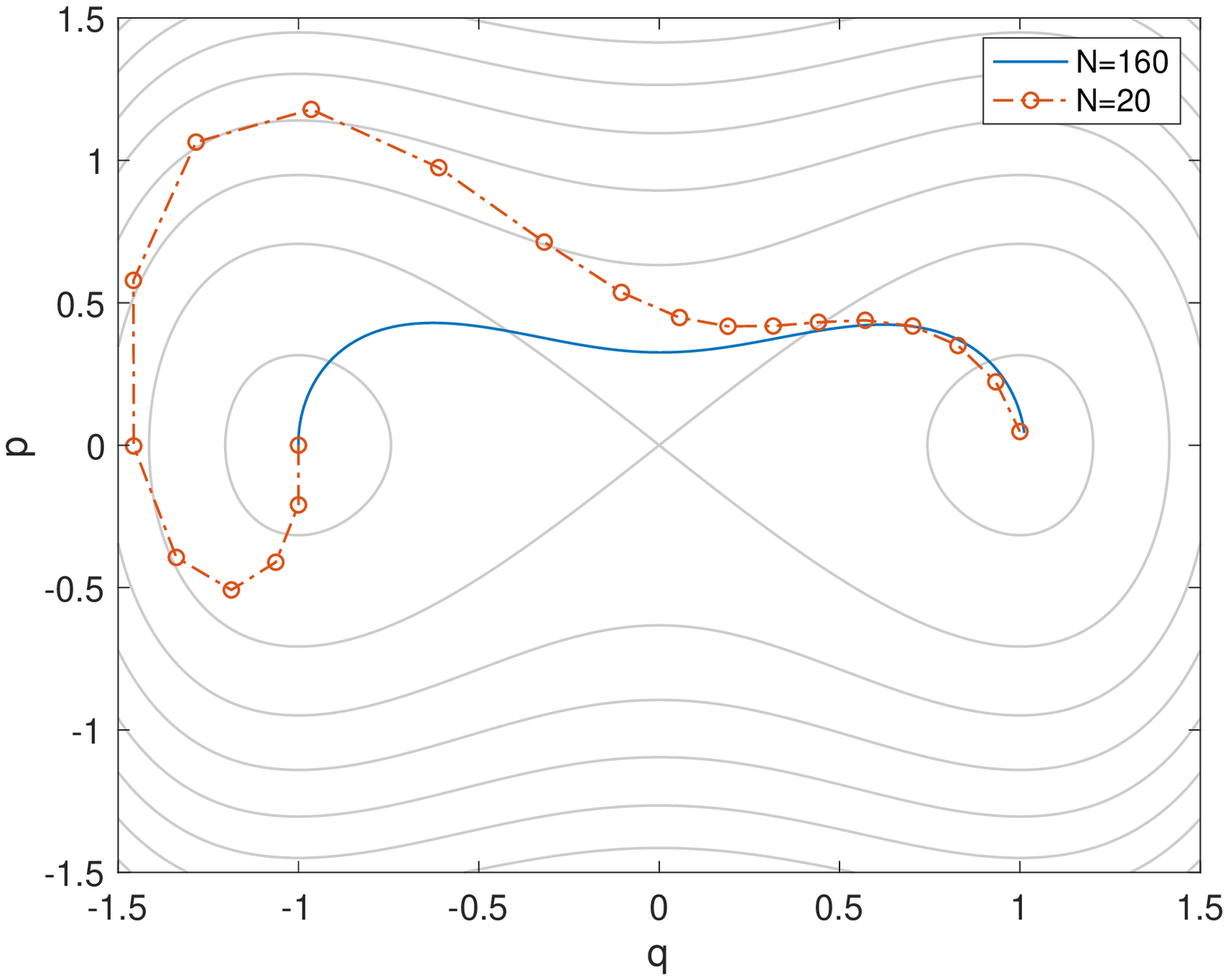}
\includegraphics[width=0.49\textwidth]{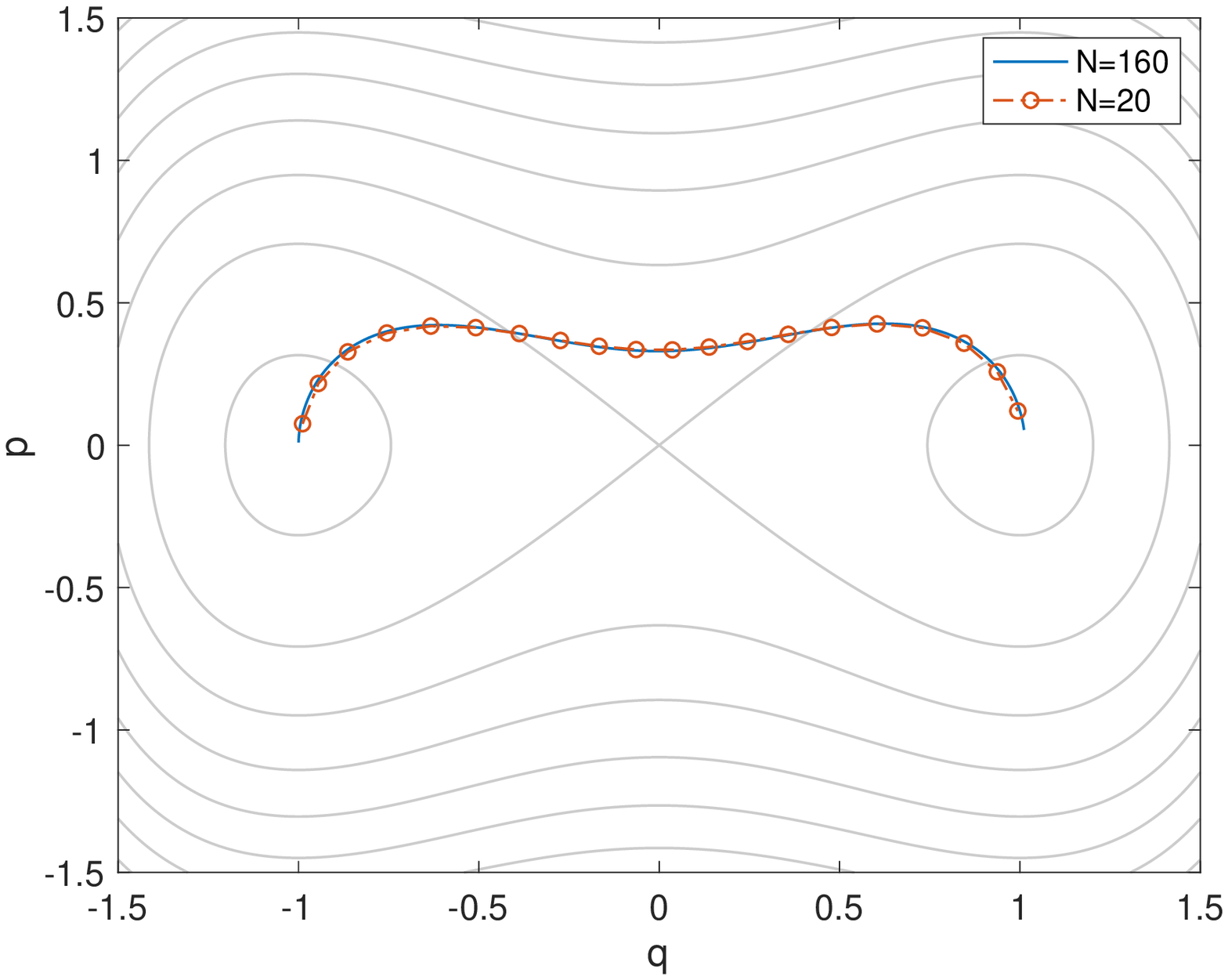}
\caption{Optimal motion in $q$--$p$ plane, computed with the symplectic Euler method (left) and implicit midpoint method (right), for $N=160$ (solid blue line) and $N=20$ (dash-dot red line). \label{fig:soln}}
\end{center}
\end{figure}

We solve the optimal control problem using the discrete regularized forward-backward sweep iteration \eqref{dMSA_x}--\eqref{dMSA_u} and the symplectic Euler scheme \eqref{Euler_x}--\eqref{Euler_u}.  We iterate until the update to the control variable $u$ is less than a prescribed tolerance
\[
	\sum_{n=0}^{N-1} \| u_n^{(k)} - u_n^{(k-1)} \| < \varepsilon,
\]
where $\varepsilon = 1e^{-8}$.
The computed optimal path $x(t) = (q(t),p(t))$ is shown as a solid blue curve on the left plot of Figure \ref{fig:soln}.  The background contours are level sets of the total energy function $E = \frac{1}{2} p^2 + \frac{1}{4} q^4 - \frac{1}{2} q^2$. The optimal control must accelerate the motion of the particle to reach an energy level above the saddle point, allowing it to cross to the potential well on the right.

\begin{figure}[!htbp]
\begin{center}
\includegraphics[width=0.60\textwidth]{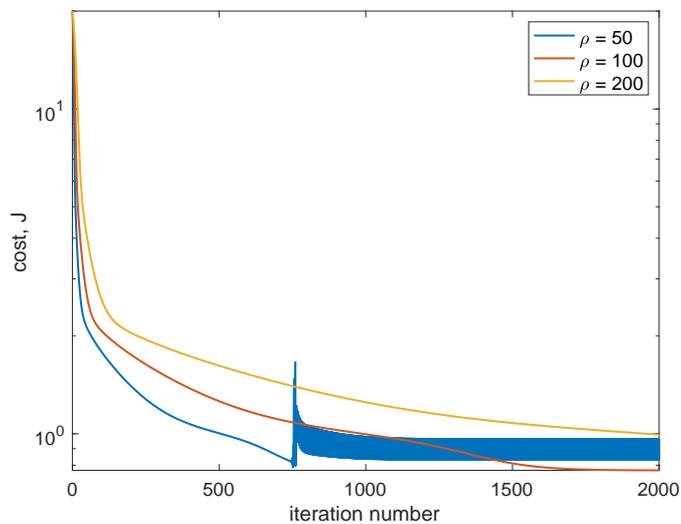}
\caption{Convergence of the cost function for the regularized forward-backward sweep iteration using the symplectic Euler method \eqref{Euler_x}--\eqref{Euler_u}, with $\rho=50$ (blue), $\rho=100$ (red) and $\rho=200$ (yellow). \label{fig:conv}}
\end{center}
\end{figure}

For this computation we chose $\rho=100$ for the regularization parameter. Convergence occurs in 4206 iterations. Figure \ref{fig:conv} shows the discrete cost function \eqref{Jtau} during the first 2000 iterations for values $\rho=50$, $\rho=100$ and $\rho=200$.
For $\rho=100$, the convergence is monotone as predicted by the theory of the previous section (cf.~\eqref{monotone}). For $\rho=50$, we observe an initial reduction in cost, which eventually oscillates and does not converge.  For $\rho=200$, the iteration converges but at a slower rate than for $\rho=100$.
Hence, our experience suggests there is a critical value of $\rho$ below which there is no convergence of the regularized forward-backward sweep iteration, and above which the convergence becomes steadily slower.

The minimal cost obtained using the symplectic Euler method and $N=160$ was $J = 0.7712$.  We also computed the optimal solution for $N=20$ time steps, shown as the red dash-dot line in the left plot of Figure \ref{fig:soln}. As noted in Section \ref{sec:VarInt}, by discretizing the Lagrangian we obtain a discrete optimal control problem for each $N$. For the case $N=20$ the optimal path deviates significantly from that for $N=160$.  Because the Lipschitz constant is larger for this solution, it was necessary to take $\rho=400$ for convergence. The optimal cost in the case $N=20$ is $J=0.7006$, which is  less than the optimal cost obtained in the case $N=160$.

We also solved the optimal control problem using the implicit midpoint rule, a second order symplectic Runge-Kutta method with $s=1$ and coefficients $a_{11} = b_1 = 1/2$.  The solutions for $N=20$ and $N=160$ are shown in the right plot of Figure \ref{fig:soln}. Here we see that the discrete optimum at low resolution is much closer to that at high resolution. The optimal costs were computed $J=0.7837$ for $N=20$ and $J=0.7769$ for $N=160$.  Both resolutions converged with $\rho=100$.

\begin{figure}[!htbp]
\begin{center}
\includegraphics[width=0.60\textwidth]{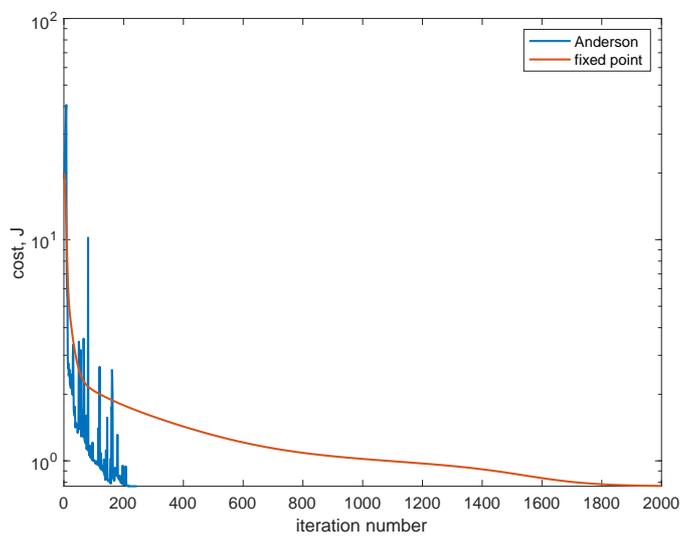}
\caption{Comparison of the Anderson accelerated (blue) and fixed point (red) iterations.  Shown are the cost functions using the symplectic Euler method \eqref{Euler_x}--\eqref{Euler_u}, with  $\rho=100$.  \label{fig:accel}}
\end{center}
\end{figure}

Although the convergence is monotone in the cost $J$ for large enough $\rho$,  the forward-backward sweep iteration may require a large number of iterations to attain a sufficiently small cost. Acceleration techniques such as Anderson acceleration \cite{WaNi11} may be employed to improve the convergence rate. We implement  \eqref{dMSA_x}--\eqref{dMSA_u} as a fixed point iteration on the control function $\B{u}$, i.e. $\B{u}^{(k+1)} = \mathcal{F} (\B{u}^{(k)})$. Subsequently we apply Anderson acceleration with restarts every three iterations.  In Figure \ref{fig:accel} we see that the cost function converges in 221 iterations (nearly a factor 20 fewer), but the cost no longer decays monotonically. See \cite{HeVa19} for a more sophisticated strategy with adaptive damping and preserving monotonicity. In our experience the choice of a good acceleration algorithm depends heavily on the problem. For instance, in other work we are investigating the use of this method for sparse control of the Cucker-Smale model with $\ell_p$ norm of the control in the running cost (see, e.g.~ \cite{bailo2018optimal}). The approach described above using Anderson acceleration works well for $p=2$, but gives no observable advantage for $p=1$.

\section{Summary\label{sec:conclusions}}
In this article we have extended the convergence proof of a regularized forward-backward sweep iteration \cite{LiChTaE18} for solving optimal control problems to the discrete setting. We showed that if the continuous problem is discretized by a symplectic partitioned Runge-Kutta pair (using a variational integrator approach), then the convergence proof of \cite{LiChTaE18} may be easily adapted. Numerical experiments with the first order, explicit symplectic Euler method and the second order implicit midpoint rule demonstrate monotonic convergence of the cost function if the regularization parameter $\rho$ is chosen large enough. For insufficiently large $\rho$ the cost undergoes bounded oscillations;  whereas for excessively large $\rho$ the convergence is slower.  In our experiments, convergence was observed even with large step sizes, however the resulting discrete optimization problem is an inaccurate approximation of the continuous problem. In an efficient implementation, the regularized forward-backward sweep iteration may be combined with an
 acceleration technique for nonlinear iterations such as Anderson acceleration \cite{WaNi11}.

\bibliographystyle{elsarticle-num}
\bibliography{Symplectic-Runge-Kutta}

\section*{Appendix}
In this appendix we prove that the bounds \eqref{discLipschitz} follow from \eqref{contLipschitz}.

Since $b_i \ge 0$, $i=1,\dots,s$,
\begin{equation}\label{dftau}
	\| f^\tau(x',u) - f^\tau(x,u) \| \le \sum_{i=1}^s b_i \| f(X_i,U_i) - f(X_i',U_i) \|,
\end{equation}
where $X_i'$ satisfies
\[
	X_i' = x' + \dt\sum_{j=1}^s a_{ij} f(X_j',U_j).
\]
Denoting $\Delta X_i = X_i - X_i'$ and using the Lipschitz condition on $f$ (cf.~\eqref{contLipschitz}), we find
\[
	\| \Delta X_i \| \le \| x - x' \| + \dt  \sum_{j=1}^s | a_{ij} | \cdot  K \| \Delta X_j \|.
\]
Denote by $|A|$ the matrix with elements $|a_{ij}|$, by $|\Delta X|$ the vector with elements $\| \Delta X_i \|$, and let $\mathbbold{1}$ be the vector of dimension $s$ with all elements equal to 1.  Then the above inequality becomes
\begin{equation}\label{stageBound1}
	(I - \dt K |A|) |\Delta X| \le \| x-x'\| \mathbbold{1}.
\end{equation}
For \emph{explicit} Runge-Kutta methods, the matrix on the left always has positive inverse given by
\[
	(I - \dt K |A|)^{-1} = \sum_{i=0}^{s-1} (\dt K |A|)^i.
\]
For \emph{implicit} Runge-Kutta methods, the matrix on the left of \eqref{stageBound1} is an M-matrix with positive inverse if we impose the step size restriction
\begin{equation}\label{steplimit}
	\dt \le (K \max_{ij} |a_{ij}|)^{-1}.
\end{equation}
In either of the above cases we find
\begin{equation}\label{stageBound}
	\| X_i - X_i' \| \le K^\tau \|x-x'\|, \qquad K^\tau =  \|  (I - \dt K |A|)^{-1} \mathbbold{1} \|_\infty.
\end{equation}
Returning to \eqref{dftau} we obtain
\[
	\| f^\tau(x',u) - f^\tau(x,u) \| \le \sum_{i=0}^s b_i K K^\dt \| x-x' \| = K K^\dt \| x-x' \|.
\]
proving the first bound in \eqref{discLipschitz}.

To prove the second bound, recall \eqref{Psi}.  Taking norms, and using the bound \eqref{contLipschitz},
\[
	\| \Psi_i \| \le 1 + \dt \sum_{j=1}^s | a_{ij} | K \| \Psi_j \|,
\]
from which we conclude that
\begin{equation}\label{Psibound}
	\| \Psi_i \| \le K^\dt.
\end{equation}
We also find
\begin{align*}
	\| \Psi_i - \Psi_i' \| &\le \dt \sum_{j=1}^s |a_{ij}| \| \prt{f}{x} (X_j,U_j) \Psi_j - \prt{f}{x} (X_j',U_j) \Psi_j' \| \\
	&= 	\dt \sum_{j=1}^s |a_{ij}| \| \prt{f}{x} (X_j,U_j) (\Psi_j - \Psi_j') + (\prt{f}{x} (X_j,U_j) - \prt{f}{x} (X_j',U_j)) \Psi_j' \| \\
	&\le \dt \sum_{j=1}^s |a_{ij}| ( K \|\Psi_j - \Psi_j' \| + K K^\dt  \| X_j - X_j' \| )\\
	&\le \dt \sum_{j=1}^s |a_{ij}| (K \| \Psi_j - \Psi_j' \| +  K(K^\dt)^2 \| x - x' \| )\\
	&\le \dt (\max_i \sum_{j=1}^s |a_{ij}|) K (K^\dt)^3 \| x - x' \|,
\end{align*}
where the last inequality follows by inverting the matrix of \eqref{stageBound1}---in the case of implicit RK methods under the step size restriction \eqref{steplimit}.
Similarly, we compute
\begin{align*}
	\| \prt{f^\dt}{x} (x,u) - \prt{f^\dt}{x}(x',u) \| &\le \sum_{i=1}^s b_i \| \prt{f}{x}(X_i,U_i) \Psi_i - \prt{f}{x}(X_i',U_i) \Psi_i' \| \\
&= \sum_{i=1}^s b_i \| \prt{f}{x}(X_i,U_i) (\Psi_i - \Psi_i') + (\prt{f}{x}(X_i,U_i) - \prt{f}{x}(X_i',U_i)) \Psi_i' \| \\
&\le \sum_{i=1}^s b_i ( K \| \Psi_i - \Psi_i' \| + K K^\dt \| X_i - X_i' \|  ) \\
& \le (\dt \max_i \sum_{j=1}^s |a_{ij}|) K^2 (K^\dt)^3 + K (K^\dt)^2)  \| x-x' \|,
\end{align*}
proving the second bound in \eqref{discLipschitz}.

The bounds on $h^\tau$ and $\prt{h^\dt}{x}$ in \eqref{discLipschitz} follow the same reasoning.

\end{document}